\documentclass[12pt]{article}

\def\beq{\begin{equation}}
\def\eeq{\end{equation}}

\def\L{\Lambda}

\def\o{\omega}
\def\O{\Omega}

\def\bc{\begin{center}}
\def\ec{\end{center}}

\def\w{\wedge}

\def\d{\partial}

\def\vuoto{\ \hfill\hbox{\vbox{\hrule\hbox{\vrule
height5pt\kern5pt\vrule height5pt}\hrule}}\par\medskip\rm}

\begin{document}

\title{Special Lagrangian Geometry in irreducible symplectic 4-folds}

\author{Alessandro Arsie \thanks{e-mail: arsie@sissa.it} \\
S.I.S.S.A. - I.S.A.S. \\
Via Beirut 4 - 34013 Trieste, Italy}
\date{}
\maketitle

\begin{abstract}
Having fixed a Kaehler class and the unique 
corresponding hyperkaehler metric, we prove that all special
 Lagrangian submanifolds of an irreducible 
symplectic 4-fold $X$ are obtained by complex submanifolds via a 
generalization of the so called Òhyperkaehler rotation trickÓ; 
thus they retain part of the rigidity of the complex submanifolds: 
indeed all special Lagrangian submanifolds of $X$ turn out to be real analytic. 
\end{abstract}

{\small MSC (1991): Primary: 53C15, Secondary: 53A40, 51P05, 53C20}

{\small Keywords: special Lagrangian submanifolds, hyperkaehler 
structures.}
\bc

{\small REF.: 75/99/FM/GEO}

\ec

\section{Introduction}

 Under the flourishing research activity on D-branes in string 
 theory, the role of special Lagrangian submanifolds in physics has 
 become more and more relevant (see for example \cite{BeBeS}) untill it 
 was eventually conjectured in 
 \cite{SYZ} that they can be considered as the cornerstones of the mirror 
 phenomenon. Indeed, D-branes are special Lagrangian submanifolds 
 equipped with a flat $U(1)$ line bundle. In physical literature, 
 special Lagrangian submanifolds of the compactification space are 
 related to physical states which retain part of the vacuum 
 supersymmetry: for this reason they are often called supersymmetric 
 cycles or BPS states.
 
 Despite their importance, there are very few  explicit examples of 
 special Lagrangian submanifolds, especially in Calabi-Yau 3-folds. 
 However, in an irreducible symplectic 4-fold (realized as a
 hyperkaehler manifold) we have a complete control of the special 
 Lagrangian geometry of its submanifolds, via a sort of Òhyperkaehler 
 trickÓ; moreover this enables us to prove that special Lagrangian 
 submanifolds retain part of the rigidity of complex submanifolds. 
 
 We first recall the following:
 
 {\bf Definition 1.1}: {\em A complex manifold $X$ is called 
 irreducible symplectic if it satisfies the following three 
 conditions:
 
 1) $X$ is compact and Kaehler;
 
 2) $X$ is simply connected;
 
 3) $H^{0}(X,\O^{2}_{X})$ is spanned by an everywhere non-degenerate 
 2-form $\o$.}
 
 In particular, irreducible symplectic manifolds are special cases of 
 Calabi-Yau manifolds (the top holomorphic form which trivializes the 
 canonical line bundle is given by a suitable power of the 
 holomorphic 2-form $\o$). In dimension 2, K3 surfaces are the only 
 irreducible symplectic manifolds, and indeed irreducible symplectic 
 manifolds appear as higher-dimensional analogues  of K3 surfaces, 
 as strongly suggested in \cite{Hu}. Unfortunately, up to now there are very few explicit 
 examples of irreducible symplectic manifolds. Indeed almost all known 
 examples turn out to be birational to two standard series of examples: Hilbert 
 schemes of points on K3 surfaces and generalized Kummer varieties 
 (both series were first studied in \cite{Beau}), but quite recently O'Grady 
 has constructed irreducible symplectic manifolds which are not 
 birational to any of the elements of the two groups (see 
 \cite{Ogra}).

Finally, let us recall from \cite{HL} the following:

{\bf Definition 1.2}: {\em Let $X$ be a Calabi-Yau n-fold, with 
Kaehler form $\o$ and holomorphic nowhere vanishing n-form $\O$. A 
(real) n-dimensional submanifold $j: \L \hookrightarrow X$ of $X$ is 
called special Lagrangian if the following two conditions are 
satisfied:

1) $\L$ is Lagrangian with respect to $\o$, i.e. $j^{*}\o=0$;

2) there exists a multiple $\O '$ of $\O$ such that $j^{*}{\rm Im}(\O 
')$=0;
one can prove (see \cite{HL}) that both conditions are equivalent to:

1') $j^{*}{\rm Re}(\O ')=Vol_{g}(\L)$.}

The condition $1')$ in the previous definition means that the 
real part of $\O '$ restricts to the volume form of $\L$, induced by 
the Calabi-Yau Riemannian metric $g$. In this way special 
Lagrangian submanifolds are considered as a type of calibrated 
submanifolds (see \cite{HL} for further details on this point).

\section{Characterization of special Lagrangian submanifolds}

In this section we will describe all special Lagrangian submanifolds 
of an irreducible symplectic 4-fold $X$  (having fixed a Kaehler 
class $[\o]$ in the Kaehler cone). The key result is the following:
 
 {\bf Theorem 2.1}: {\em Every connected special Lagrangian submanifold of an 
 irreducible symplectic 4-fold is also bi-Lagrangian, in the sense 
 that it is Lagrangian with respect to two different symplectic 
 structures.} 
 
 {\bf Proof}: Let us a fix a Kaehler class on the irreducible 
 symplectic 4-fold $X$. By Yau's Theorem this determines a unique 
 hyperkaehler metric $g$. Choose a hyperkaehler structure 
 $(I,J,K)$ compatible with the metric $g$ (notice that the triple 
 $(I,J,K)$ is not uniquely determined) and consider the associated 
 symplectic structures $\o_{I}(.,.):=g(I.,.)$, $\o_{J}(.,.):=g(J.,.)$ 
 and $\o_{K}(.,.):=g(K.,.)$.
 
  Consider a special Lagrangian 
 submanifold $\L$ in the complex structure $K$ (this is not 
 restrictive, since $(I,J,K)$ is not uniquely determined); that is assume that 
 $\L$ is calibrated by the real part of the holomorphic (in the 
 structure $K$) 4-form:
 \beq
 \label{4-form}
 \O_{K}:=\frac{1}{2!}(\o_{I}+i\o_{J})^{2}.
 \eeq
 Notice that the real and immaginary part of $\O_{K}$ are then given 
 by:
 
 \beq
 \label{reim}
{\rm  Re}(\O_{K})=\frac{1}{2}(\o^{2}_{I}-\o^{2}_{J}) \quad 
 {\rm Im}(\O_{K})=\o_{I}\w\o_{J}.
 \eeq
 
 Obviously, by the property of being special Lagrangian we have 
 that $\L$ is Lagrangian with respect to $\o_{K}$. We will prove that 
 having fixed the calibration, if $\L$ is not Lagrangian also with 
 respect to $\o_{I}$, then it is necessarily Lagrangian with 
 respect to $\o_{J}$. 
 First we work locally and consider $V:=T_{p}\L$ ($p\in \L$), spanned 
 by $(w_{1},w_{2},w_{3},w_{4})$. Since $\L$ is assumed not to be Lagrangian with 
 respect to $\o_{I}$, we have to deal with two cases.
 
 First case: $V$ is a symplectic vector space for the structure $\o_{I}$. In this 
 case we can choose a symplectic basis for $V$ and this can always 
 be chosen to be of the form ${v_{1},Iv_{1},v_{2},Iv_{2}}$. Then 
 $V$ is Lagrangian in the symplectic structure $\o_{J}$; indeed 
 $\o_{J}(v_{1},Iv_{1})=g(Jv_{1},Iv_{1})=g(IJv_{1},-v_{1})=-\o_{K}(v_{1},v_{1})=0$;
 analogously for $\o_{J}(v_{2},Iv_{2})$; 
 $\o_{J}(v_{1},Iv_{2})=g(Jv_{1},Iv_{2})=-\o_{K}(v_{1},v_{2})=0$ since 
 $v_{1},v_{2}$ belong to a Lagrangian subspace of $\o_{K}$, and 
 analogously for $\o_{J}(v_{2},Iv_{1})=-\o_{K}(v_{2},v_{1})=0$. Thus 
 $V$ is also Lagrangian for the symplectic structure $\o_{J}$.
 
 Second case: $V$ is neither symplectic nor Lagrangian for the 
 structure $\o_{I}$. Notice $V$ can not be 
 symplectic with respect to $\o_{J}$, otherwise by the first case it 
 would be Lagrangian in the strucutre $\o_{I}$; moreover we can assume 
 that $V$ is not Lagrangian with respect to $\o_{J}$, otherwise there 
 is nothing to prove. So in this case $V$ is neither Lagrangian nor 
 symplectic in the structure $\o_{I}$ {\em and} in the structure $\o_{J}$. 
 This means that $V$ contains a symplectic 2-plane $\pi$ with respect to $\o_{I}$
 and a symplectic 2-plane $\rho$ with respect to $\o_{J}$. 
 Indeed, consider $v_{1}\in V$; since $V$ is not Lagrangian in the 
 structure $\o_{I}$, there exists $v_{2}\in V$ such that 
 $\o_{I}(v_{1},v_{2})\neq 0$ and this implies that the vector 
 subspace $\pi$ spanned by $(v_{1},v_{2})$ is a symplectic vector 
 space for $\o_{I}$, which can not be extended to all $V$. The same 
 reasoning applies in the structure $\o_{J}$.
 
 We prove that this can not happen, since it violates the calibration 
 condition. We have to distinguish three different subcases according 
 to the intersection of $\pi$ with $\rho$.
 
 First subcase: $\pi$ and $\rho$ have zero intersection.
 If this happens we can always choose a basis of $V$ 
 of the form $(v_{1}, Iv_{1}, v_{2}, Jv_{2})$. Write $\pi$ for 
 the 2-plane spanned by ${v_{1}, Iv_{1}}$ and $\rho$ for that spanned 
 by $v_{2}, Jv_{2}$, so that $V=\pi\oplus \rho$. Indeed, since $V$ is 
 not Lagrangian with respect to $\o_{I}$, it has to contain a 
 symplectic 2-plane like $\pi$, and similarly for $\rho$ and $\o_{J}$. 
 Moreover, since $V$ is not symplectic with respect to $\o_{I}$, it 
 turns out that the symplectic 2-plane $\pi$ can not be completed to 
 a symplectic basis of $V$, so that $V$ has to contain an isotropic 
 2-plane for $\o_{I}$, which is $\rho$. The same reasoning (with the 
 roles reversed) applies obviously to the symplectic structure $\o_{J}$.
 Hence, in this case we have:
 \[ 
 2{\rm Re}(\O_{K})|_{V}=(\o^{2}_{I}-\o^{2}_{J})(v_{1},Iv_{1},v_{2},Jv_{2})=
\o_{I}(v_{1},Iv_{1})\o_{I}(v_{2},Jv_{2})-\]
\[
\o_{I}(v_{1},v_{2})\o_{I}(Iv_{1},Jv_{2})+\o_{I}(v_{1},Jv_{2})\o_{I}(Iv_{1},v_{2})-
\o_{J}(v_{1},Iv_{1})\o_{J}(v_{2},Jv_{2})+\]
\[
\o_{J}(v_{1},v_{2})\o_{J}(Iv_{1},Jv_{2})
-\o_{J}(v_{1},Jv_{2})\o_{J}(Iv_{2},v_{2})=0,\]
using the defining relations of $\o_{I},\o_{J},\o_{K}$, the 
quaternionic relation $IJ=K$, the invariance of $g$ and the fact that $V$ is Lagrangian with 
respect to $\o_{K}$. So this subcase is not consistent with the 
calibration property. 

Second subcase: $\pi$ and $\rho$ have a 1-dimensional intersection 
spanned by a vector $v_{1}$. In this case we can choose a basis of 
$V$ of the form $(v_{1},Iv_{1},Jv_{1},w)$ ($\pi$ is spanned by 
$(v_{1},Iv_{1})$, while $\rho$ is spanned by $(v_{1},Jv_{1})$). Again 
by the same computation of the previous subcase one shows that this 
configuration is not compatible with the calibration.

Third subcase: Finally $\pi$=$\rho$ can not clearly happen, since 
otherwise one can choose a basis of $\pi$ equal to $(v_{1},Iv_{1})$, but then, 
in this basis $\o_{J}$ is identically vanishing, contrary to the 
assumption that $\rho=\pi$ is a symplectic 2-plane also for $\o_{J}$. 

Since the second case can never happen $V$ has to be Lagrangian also 
with respect to $\o_{J}$.

Up to now, we have worked only locally; to conclude the proof it is 
necessary to show that if $T_{p}\L$ is Lagrangian with respect to 
$\o_{J}$, then it can not be possible that $T_{q}\L$ is Lagrangian 
with respect to $\o_{I}$, for a different $q\in \L$. Notice that any 
tangent space to $\L$ can not be Lagrangian with respect to {\em 
both} $\o_{I}$ and $\o_{J}$, otherwise it would violates the 
calibration condition. Consider now the following smooth sections of 
$\bigwedge^{2} T^{*}\L$:
\[
\begin{array}{cccc}
          \alpha_{I,J}: & \L & \rightarrow & \bigwedge^{2} T^{*}\L \\ 
                &  p & \mapsto &\o_{I,J}{|T_{p}\L} \\
\end{array}
\]
and the zero section $s_{0}: \L \rightarrow \bigwedge^{2} T^{*}\L$.     
Obviously, $s_{0}(\L)$ is closed in $\bigwedge^{2} T^{*}\L$, and by 
the previous reasoning $\L$ can be decomposed as 
$\L=\alpha_{I}^{-1}(s_{0}(\L))\cup \alpha_{J}^{-1}(s_{0}(\L))$, that 
is as the disjoint union of two proper closed subsets. But this is clearly 
impossible, since $\L$ is connected, and this implies that one of the 
two closed subset is empty, so $\L$ is bi-Lagrangian. \vuoto
 
The previous theorem is important in view of the following:

{\bf Corollary 2.1}: {\em Every (connected, compact and without border) 
special Lagrangian submanifold $\L$  of a 
hyperkaehler 4-fold $X$ can be realized as a complex submanifold, 
via hyperkaehler rotation of the complex structure of $X$.}

{\bf Proof}: Let $\L$ be a special Lagrangian submanifold of $X$ in 
the complex structure $K$. Then by definition 
${\rm Re}(\O_{K})_{|\L}={\rm Vol}_{g}(\L)$, but by the previous theorem, since 
$\o_{J}|_{\L}=0$ this means:
\beq
\label{vol} 
 {\rm Vol}_{g}(\L)=\frac{1}{2}\int_{\L}\o^{2}_{I}.
 \eeq
By Wirtinger's theorem, since $\L$ is assumed to be compact and 
without border, condition (\ref{vol}) is equivalent to say that $\L$ is a 
complex submanifold of $X$, in the complex structure $I$, that is 
performing a hyperkaehler rotation. Notice that in the complex structure 
$I$, $\L$ is still a Lagrangian submanifold with respect to $\o_{K}$ 
and $\o_{I}$, so it is Lagrangian with respect to the holomorphic (in 
the structure $I$) 2-form $\O_{I}:=\o_{J}+i\o_{K}$. \vuoto 

Collecting the results so far proved, we can show that special 
Lagrangian submanifolds of $X$ are particularly ÒrigidÓ:

{\bf Proposition 2.1}: {\em Any (connected, compact and without border) special
 Lagrangian submanifold $\L$  of 
a hyperkaehler 4-fold $X$ is real analytic.}

{\bf Proof}: Let $\L$ be a special Lagrangian submanifold of $X$, 
having fixed some complex structure on $X$, let us say $K$; then, by Corollary 2.1 there 
exists a new complex structure, let us say $I$, in which $\L$ is 
holomorphic, that is, it is locally given by:
\[f_{1}(z_{1},\ldots,z_{4})=0 \quad {\rm and} \quad 
    f_{2}(z_{1},\ldots,z_{4})=0. \]  
Now observe that coming back to the original complex structure $K$, we 
induce an {\em analytic} change of coordinates from the holomorphic 
coordinates $z^{i}$ ($I\frac{\d}{\d z^{i}}=i \frac{\d}{\d z^{i}}$) to 
new holomorphic coordinates $w^{i}$ ($K\frac{\d}{\d w^{i}}=i 
\frac{\d}{\d w^{i}}$) such that locally:
\beq 
\label{analy}
z^{i}=c_{1}w^{i}+c_{2}\bar{w}^{i} \quad 
\bar{z}^{i}=d_{1}w^{i}+d_{2}\bar{w}^{i}, 
\eeq 
for some complex constants $c_{j},d_{j}$. Thus in the complex 
structure $K$ the special Lagrangian submanifold $\L$ is given 
by  $f_{j}(c_{1}w^{i}+c_{2}\bar{w}^{i},d_{1}w^{i}+d_{2}\bar{w}^{i})=0$
which is again the zero locus of a set of functions analytic in 
$w^{i},\bar{w}^{i}$. \vuoto

Quite naturally, the action of the hyperkaehler rotation can be 
extended also to the holomorphic functions defined on complex 
submanifolds $S$ of $X$; in particular we have an action of the 
hyperkaehler rotation on the structure sheaf ${\cal O}_{S}$ (here, as 
always, we identify ${\cal O}_{S}$ with its direct image $j_{*}{\cal O}_{S}$, 
where $j:S\hookrightarrow X$ is the holomorphic embedding).
We are thus led to give the following:

{\bf Definition 2.2}: {\em Let $\L$ be a special Lagrangian 
submanifold of a hyperkaehler 4-fold $X$ (in the complex structure $K$). Then we define the special 
Lagrangian structure sheaf ${\cal L}_{\L}$ as the sheaf obtained 
by the action of the hyperkaehler rotation on the structure sheaf 
${\cal O}_{\L}$ of $\L$, as a complex Lagrangian submanifold of $X$, 
(in the structure $I$).}

\section{Concluding remarks}

It is important to remark that all previous results are true also for 
special Lagrangian submanifolds of K3 surfaces, but their proof is 
completely trivial in that case.

Another observation is related to {\em singular} Lagrangian 
submanifolds: indeed, by the previous results, it turns out that we can 
also give examples of special Lagrangian {\em subvarieties}, obtained 
via hyperkaehler rotation of Lagrangian complex subvarieties. On the 
other hand, contrary to the case of the corresponding submanifolds, 
 we can not expect that {\em all} special Lagrangian subvarieties are 
 obtained in this way, and consequently we can not expect that all 
 special Lagrangian subvarieties are real analytic. Indeed, there are 
 examples (compare \cite{HL}) of singular special Lagrangian 
 submanifold in $C^{n}$ which are only smooth, but not real 
 analytic. 
 
 The discussion about singular Lagrangian submanifolds leads us to 
 comment on the mirror symmetry construction suggested in \cite{SYZ}. 
  Indeed, according to the recipe of \cite{SYZ},  
  any  Calabi-Yau $X$, admitting a mirror $\hat{X}$, has a peculiar
  fibre space structure: on a physical ground it is argued that $X$ can be 
 realized as the total space of a fibration in special Lagrangian 
 tori. Unfortunately, there are very few examples of such realization: 
 in particular, as far as we know, there is only one (partial) example for 
 Calabi-Yau 3-folds of the so called Borcea-Voisin type (see \cite{GW}). 
 Instead, in the case of irreducible symplectic 
 projective manifolds the situation is completely different. Indeed, 
 a recent result of Matsushita (see \cite{Mat1} and \cite{Mat2}) 
 shows that for any fibre space structure $f:X\rightarrow B$ of a 
 projective irreducible symplectic manifold $X$, with projective base 
 $B$, the generic fibre $f^{-1}(b)$ is an Abelian variety (up to 
 finite unramified cover), and it is also Lagrangian with respect to 
 the non degenerate holomorphic 2-form $\O$; moreover, in the case of 
 4-folds one can prove that the generic fibre is an Abelian surface 
 and $f$ is equidimensional, (i.e. all irreducible components of the 
 fibres have the same dimension). By Corollary 2.1 it turns out that 
 this fibre space structure can also be realized as a special 
 Lagrangian torus fibration; moreover, in this case all special 
 Lagrangian fibres, even the singular ones, are analytic, since they 
 are obtained by performing a hyperkaehler rotation starting from Lagrangian 
 Abelian surfaces. So, in these cases, we have special Lagrangian 
 torus fibration in which all fibres are analytic: one 
 can hope to understand the degeneration types of singular special 
 Lagrangian tori, moving from these constructions.
  
 Explicit examples of projective irreducible symplectic 4-folds, 
 fibered over a projective base have been constructed by Markuschevich 
 in \cite{Mar1} and \cite{Mar2}. One of this constructions is the 
 following: consider a double cover $\pi: S \rightarrow 
 P^{2}$ of the projective plane, ramified along a smooth sextic 
 $C \hookrightarrow P^{2}$ ($S$ is then realized as a K3 
 surface). Since any line in $P^{2}$ will intersect generically the sextic $C$ in six distinct 
 point, we have that the covering $\pi:  S \rightarrow P^{2}$ determines a (flat) family of hyperelliptic curves over the 
 dual projective plane $f: {\cal X} \rightarrow P^{2}$. Then 
 the Altmann-Kleiman compactification of the relative Jacobian of the 
 family turns out to be a simplectic  projective irreducible 4-folds, 
 fibered over $P^{2}$, and in fact all fibres are Lagrangian 
 Abelian varieties.

Finally, we believe that our characterization of special Lagrangian 
submanifolds of irreducible symplectic 4-folds can be extended also to 
higher dimensional irreducible symplectic manifolds: to this aim 
notice that the proof we have given becomes longer and longer, since 
one has to deal with new cases and subcases. It would be nice, instead, to find 
out a sort of inductive argument, which works for all dimensions.

\end{document}